\documentclass[english]{amsart}

\usepackage{babel}
\usepackage{amstext}
\usepackage{amsmath}
\usepackage{amsfonts}
\usepackage{latexsym}
\usepackage{ifthen}
\usepackage{xypic}
\xyoption{all}
\newcommand{\sO}{{\mathcal O}}
\newcommand{\PN}{{\mathbb P}}

\newcommand{\KQ}{{\mathbb Q}}
\newcommand{\Pic}{{\rm Pic}}

\newcommand{\lra}{\longrightarrow}
\newcommand{\KC}{{\mathbb C}}

\newcommand{\X}{\mathcal X}
\newcommand{\KZ}{\mathbb Z}

\newcommand{\W}{\mathcal W}
\newcommand{\Y}{\mathcal Y}

\newtheorem{lemma1}[equation]{}

\newenvironment{example}{\begin{lemma1}{\bf Example.}\rm}{\end{lemma1}}
\newenvironment{theorem}{\begin{lemma1}{\bf Theorem.}}{\end{lemma1}}

\newenvironment{proposition}{\begin{lemma1}{\bf Proposition.}}{\end{lemma1}}
\newenvironment{corollary}{\begin{lemma1}{\bf Corollary.}}{\end{lemma1}}

\begin{document}

\title{Terminal Fano threefolds and their smoothings}
\author{Priska Jahnke}
\author{Ivo Radloff}
\address{Mathematisches Institut \\ Universit\"at Bayreuth \\ D-95440 Bayreuth/Germany}
\email{priska.jahnke@uni-bayreuth.de}
\email{ivo.radloff@uni-bayreuth.de}
\thanks{The authors gratefully acknowledge support by the Schwerpunkt program {\em Globale Methoden in der komplexen Geometrie} of the Deutsche Forschungsgemeinschaft.}
\date{\today}
\maketitle

\section*{Introduction}

Let $X$ be some irreducible projective variety. By a {\em smoothing of $X$} we
mean a flat projective morphism
  \[\pi: \X \lra \Delta\]
onto the unit disc $\Delta$, such that $\X_0 \simeq X$ while $\X_t$ is
smooth for $t \not=
0$. The total space $\X$ is a reduced
complex space but not necessarily smooth. We are interested in
smoothings in the case of singular (almost) Fano
threefolds $X$. Recall that $X$ is called Fano if some multiple of $-K_X$ is
free and ample (big and nef in the case almost Fano). 

The existence of a smoothing is known here in the following cases:
 \begin{enumerate}
   \item $X$ is a Gorenstein Fano threefold with at most terminal
     singularities (\cite{Nam}). Here $\X_t$ is Fano.
   \item $X$ is an almost Fano threefold with at most $\KQ$--factorial
     terminal Gorenstein singularities (\cite{Mi}). Here $\X_t$ is almost Fano.
\end{enumerate}
Associated to any (almost) Fano threefold we have several discrete
quantities, among them the degree, by definition the top intersection number $(-K_X)^3 > 0$,
and the Picard number
  \[\rho(X) = {\rm rk}_{\KZ}\Pic(X).\]
The degree is of course constant in a smoothing. The
Picard number may jump. 

We show that the Picard number does not jump in the above
two cases and give some examples around smoothings of
Gorenstein Fanos in the case of canonical singularities. 

\section{Terminal Gorenstein Singularities}
\setcounter{equation}{0}

Throughout this section, $\X$ is a smoothing of a terminal Gorenstein
(almost) Fano threefold $X = \X_0$. Simplest examples are terminal
degenerations of hypersurfaces of degree at most four in $\PN_4$.

\begin{proposition} \label{factorial}
  $\X$ is normal having at most isolated terminal factorial Gorenstein singularities.
\end{proposition}

\begin{proof}
  Terminal Gorenstein singularities in dimension $3$ are isolated hypersurface
  singularities (\cite{YPG}). Through any point $P \in \X$ we have the
  reduced Cartier divisor $\X_t$ having at most singularities of this type. Sequence
  \[N^*_{\X_t/\X} \lra \Omega^1_{\X}|_{\X_t} \lra \Omega^1_{\X_t} \lra
  0\]
  shows $e_P(\X) = 4$ or $5$ for the local embedding dimension. Either
  $P$ is a smooth point or an analytic isolated hypersurface
  singularity. We conclude that $\X$ is Cohen--Macaulay, normal and
  Gorenstein. By Inversion of Adjunction (\cite{KoMo},
Theorem 5.50) $(X, 0)$ klt implies $(\X, X)$ plt. The singularities
of $\X$ are hence terminal. Finally, Grothendieck's proof of Samuel's Conjecture (\cite{Gro},
Corollaire 3.14) implies that all stalks are
factorial, i.e., ${\mathcal C}l(\sO_{\X, P}) = 0$ (compare \cite{Rem},
p.85 for the analytic case). Then $\X$ is factorial.
\end{proof}

\

Recall the definition of a {\em simultaneous resolution}. It means a diagram
   \[\xymatrix{\X' \ar[rr]\ar[rd]_{\pi'} & & \X \ar[ld]^{\pi} \\ &
     \Delta &}\]
where each $\X'_t$ is smooth and $\X'_t \to \X_t$ is birational for
all $t$. The next Corollary follows from the fact that $\X$ is factorial:

\begin{corollary}
  A simultaneous resolution does not exist (unless $X$ is smooth). 
\end{corollary}

We come to the Picard groups. The Kawamata--Viehweg vanishing Theorem says $H^i(\X_t, \sO_{\X_t}) =
0$ for $i > 0$ and all $t$. By the Leray spectral sequence $H^i(\X,
\sO_{\X}) = 0$ for $i > 0$. The exponential sequence gives
  \[\Pic(X) \simeq H^2(X, \KZ) \quad \mbox{and} \quad \Pic(\X) \simeq
  H^2(\X, \KZ).\]
The inclusion $X \hookrightarrow \X$ is a deformation retract. Hence
  \begin{equation} \label{first}
  \Pic(X) \simeq H^2(X, \KZ) \simeq H^2(\X, \KZ) \simeq \Pic(\X)
 \end{equation}
by restriction and these groups are torion free. A priori the Picard
number of $\X_t$ might be different from the Picard number of $X$.

\begin{theorem} \label{pic}
  We have
   \[\Pic(X) \simeq \Pic(\X) \quad \mbox{and} \quad \Pic(\X_t) \simeq \Pic(\X).\]
\end{theorem}

In the case $X$ smooth this just follows from (\ref{first}) and the
fact that the smoothing is differentiably trivial. 
Example~\ref{jump} below illustrates some steps in the following proof.
\begin{proof}
Only the second isomorphism has to be proved. The vanishing $H^{2,0}(\X_t) = 0$ implies that the primitive part of $H^2(\X_t)$ defines a rational polarized variation of Hodge structure
of pure type $(1, 1)$. From the fact that the isometries of an
integral lattice form a finite group, we get that the corresponding monodromy is
finite, i.e., the image of
    \[\rho: \pi_1(\Delta^*) \simeq \KZ \lra H^2(\X_t, \KZ)\]
is a finite cyclic group of some order $N$. After base change
 \[\xymatrix{\X' \ar[rr]\ar[d] && \X \ar[d] \\
    \Delta \ar[rr]^{s \mapsto s^N = t} && \Delta}\]
the monodromy becomes trivial. We have $\X'_0  \simeq X$, i.e., $\X'$ is another smoothing of $X$. By
  Proposition~\ref{factorial}, $\X'$ is normal with at most
  factorial terminal Gorenstein singularities. Our aim is to show 
  \begin{equation} \label{toshow}
      \Pic(\X')\otimes \KQ \simeq H^2(\X', \KQ) \simeq H^2(\X'_t,
      \KQ) \simeq \Pic(\X'_t) \otimes \KQ.
    \end{equation}
  Then, as $H^2(\X, \KQ) \simeq H^2(X,
  \KQ) \simeq H^2(\X', \KQ)$, we find $H^2(\X, \KQ) \simeq H^2(\X_t, \KQ)$. The image of
$H^2(\X, \KQ) \lra H^2(\X_t,\KQ)$ consists of monodromy invariant
classes. If this map is an isomorphism, the monodromy must be trivial.
Then (\ref{toshow}) shows the second identity for rational line bundles,
i.e., $\Pic_{\KQ}(\X) \simeq \Pic_{\KQ}(\X_t)$.

\vspace{0.2cm}

For simplicity denote $\X'$ by $\X$. In order to prove (\ref{toshow})
we have to show that any rational cohomology class $\xi_t \in H^2(\X_t, \KQ)$
comes by restriction from a class in $H^2(\X, \KQ)$. The monodromy being trivial, we may think of
$\xi_t$ as being defined on $\X^* = \X \backslash X$, varying holomorphically
with $t$. Notice that even if some multiple of $\xi_t$ is represented by an effective
divisor, we cannot just take the closure in $\X$, in the analytic
setting this closure need not be a
divisor. We will use the Local Invariant Cycle Theorem (LICT) instead 
to extend $\xi_t$ to $\X_{reg}$, and then make use of the fact that $\X$ is
factorial to extend it to $\X$.

\vspace{0.2cm}

We first recall the statement of the LICT. Let $\pi: \W \lra \Delta$
be a semistable degeneration of K\"ahler manifolds. As above, any
rational monodromy invariant cohomlogy class on some $\W_t$, $t \not= 0$,
extends to $\W^*$ by the degeneration of the Leray spectral sequence. The LICT says that it even 
extends to the whole of $\W$, i.e., we
have a surjection:
 \[H^k(\W, \KQ) \lra H^k(\W_t, \KQ)^{inv} \lra 0.\]
The LICT is due to Clemens in the K\"ahler case (\cite{C}. See also \cite{D},
\cite{V}. Semistable degeneration might not be necessary here,
compare \cite{C}, p.230).

\vspace{0.2cm}

In order to apply the LICT let $\Y$ be a log resolution of $(\X, X)$ and $\W$ be a semistable
reduction of $\Y$
\[\xymatrix{\W \ar@{..>}[r]\ar[rd] & \Y_b \ar[d]\ar[r] & \Y \ar[d] \\
                              & \Delta \ar[r]^{t \mapsto t^b} & \Delta}\]
The pullback of $\xi_t$ to $\W^*$ comes by
restriction from a class $\xi_{\W} \in H^2(\W, \KQ)$. The rational map from $\W$
to $\X$ is a Galois covering over $\X_{reg}$. The Galois trace of
$\xi_{\W}$ induces a class $\xi \in H^2(\X_{reg}, \KQ)$ and extends $\xi_t$ to $\X_{reg}$.

As $\X$ is Cohen--Macaulay, $H^2(\X_{reg}, \sO_{\X_{reg}}) = 0$. This
follows from the long exact sequence of cohomology with support. The
exponential sequence shows once again that a
multiple of $\xi$ corresponds to a line bundle on
$\X_{reg}$. As $\X$ is normal and factorial, this line bundle comes
from a line bundle on $\X$. Its Chern class gives the desired
extension of $\xi_t$ to $\X$ and shows (\ref{toshow}).

\vspace{0.2cm}

Now $\X' = \X$ by (\ref{toshow}) and $\Pic(\X)
\hookrightarrow \Pic(\X_t)$ is an isomorphism after tensorizing with
$\KQ$. We now prove surjectivity over $\KZ$. Notice that the above
argumentation involves a trace not defined over $\KZ$.

Let $L_t \in
\Pic(\X_t)$. After tensorizing with $\sO_{\X}(-mK_{\X})$ for some $m
\gg 0$ we may assume that each $L_t$ is effective and consider a
family $D_t \in
|L_t|$. Again the closure of $D_t$ in $\X$ might not be a divisor. But
take a pointwise limit $D_{t_n}$ as $t_n \to 0$. The limit
$D_0 \in Cl(X)$ depends on $(t_n)$.

By (\ref{toshow}), $L^{\otimes r}_t$ comes by restriction of a line
bundle $B$ on $\X$. Then $rD_0$ is Cartier. As $X$ is Gorenstein
with terminal singularities, any $\KQ$--Cartier Weil divisor is
Cartier (\cite{Kaw}, Lemma~5.1). By (\ref{first}) $D_0$ induces a
line bundle $L \in \Pic(\X)$. As $L^{\otimes r} = B$ we get $L|_{\X_t}
\simeq L_t$, i.e., $L$ is the desired extension of $L_t$ to $\X$.
\end{proof}

\section{Canonical Gorenstein Singularities}
\setcounter{equation}{0}

In this section we collect examples. From now on $X$ is an (almost)
Fano variety with at most canonical Gorenstein singularities. Here a smoothing need not exist (Example~\ref{nonex}). If it exists, the Picard number may jump
(Example~\ref{jump}). Different smoothing are possible (Example~\ref{2smoothings}). 

Recall the definition of a cDV (compound DuVal) point. If $P \in X$ is a
threedimensional canonical singularity, then a general hyperplane section
through $P$ either has a DuVal or an elliptic singularity. In the
first case, $P$ is called cDV, in the latter
non--cDV. In Example~\ref{nonex} and Example~\ref{2smoothings} we have
non-cDV singularities.

\begin{example} \label{nonex}
The system of cubics embeds $\PN_2$ into $\PN_9$. Let $X$ be the cone
in $\PN_{10}$. The vertex is a non--cDV canonical point. Its blowup
yields the almost Fano threefold $Y = \PN(\sO_{\PN_2} \oplus \sO_{\PN_2}(3))$ of degree
$(-K_X)^3 = ({-}K_Y)^3 = 72$ (see \cite{JPR},
table A.3, no.1). A smooth Fano
threefold of degree $72$ does not exist, implying that $X$ is not smoothable. 
\end{example}

Examples of smoothings of canonical Fanos where the
Picard number jumps are easy. A simple one in the case of
surfaces is the degeneration of $\PN_1
\times \PN_1$ into the quadric cone in $\PN_3$; for a threefold take the product
with $\PN_1$. This example and the simultaneous resolution
quite nicely illustrates some steps in the proof of
Theorem~\ref{pic}:

\begin{example} \label{jump}
 Consider the following concrete degeneration of $\PN_1 \times \PN_1$ into the quadratic cone $Q_0$
    \[\X = \{tw^2 + x^2 + y^2 + z^2 = 0\} \subset \PN_3 \times \Delta.\]
The total space $\X$ is smooth but the monodromy is non--trivial: choose a branch
of the square root. For fixed $t \not= 0$, where this branch is defined, the two rulings $f_1, f_2$ of $\X_t \simeq \PN_1 \times \PN_1$ are given by the lines
    \[\pm \sqrt{t}\;w = \sqrt{{-}1}\;x \quad \mbox{and} \quad y = \sqrt{{-}1}\;z.\]
Going around the origin, $\sqrt{t}$ becomes $-\sqrt{t}$, i.e., the
  monodromy interchanges the rulings:
   \[T_2 = \left(\begin{array}{cc}
       0 & 1 \\
       1 & 0
   \end{array}\right)\quad \mbox{on } H^2(\X_t, \KZ) = \KZ[f_1] \oplus \KZ[f_2].\]
In order to trivialize the monodromy we make a base change  $t \mapsto t^2$. We get $\X'$ given by
   \[t^2w^2 + x^2 + y^2 + z^2 = 0 \quad \mbox{in } \PN_3 \times \Delta.\]
  Now the monodromy is trivial. The central fiber is still $\simeq
  Q_0$. But $\X'$ is no longer smooth, it has
  one terminal point in the central fiber at $[1
  :0:0:0]$, lying over the vertex of the cone $Q_0$. This singularity is {\em not factorial},
  and we cannot continue as above.

  The problem is the now globally on $\X'$ defined (Weil-)
  divisor
    \[tw = \sqrt{{-}1}\;x \quad \mbox{and} \quad y = \sqrt{{-}1}\;z\]
  which is not $\KQ$--Cartier. The blowup of this divisor gives a small
  resolution $\Y \to \X'$ (a simultaneous resolution of the canonical family
  $\X$).

The central fiber of $\Y$ is the second Hirzebruch surface $\Sigma_2$, i.e., $\Y$ comes from letting degenerate the non--split extension
    \[0 \lra \sO_{\PN_1} \lra \sO_{\PN_1}(1) \oplus \sO_{\PN_1}(1) \lra \sO_{\PN_1}(2)
    \lra 0\]
to the splitting one, parametrised by $\KC = H^1(\PN_1,
\sO_{\PN_1}(-2))$. On $\Y$ we can argue as above, finding that the general
  and the central fiber have the same Picard number, which is clear.
\end{example}

\begin{example} \label{2smoothings}
  Let $X$ be the cone over the del Pezzo surface $S$ of degree
  six. This is a Gorenstein Fano threefold with canonical
  singularities and Picard number one. Here we have two smoothings, one with general
  fiber $\PN(T_{\PN_2})$ and one with fiber $\PN_1 \times \PN_1 \times
  \PN_1$. 

Indeed, we may think of $S$ as a hyperplane section in the half anticanonical system
  of $\PN(T_{\PN_2})$ or of $\PN_1 \times \PN_1 \times \PN_1$. Any
  smooth projective manifold flatly degenerates to the cone over its general hyperplane section.

With this remark it is easy to give explicit equations: the anticanonical system of $S$
embeds $S$ into $\PN_6$. Think of $S$ as the blowup $Bl_{p_1,
  p_2, p_3}(\PN_2)$. After some
projective transformation we may assume $p_1 = [1:0:0]$, $p_2 =
[0:1:0]$, $p_3 = [0:0:1]$. The anticanonical embedding is then given by all cubics in $\PN_2$ through these points:
  \[x_0 = x^2y, \; x_1 = x^2z, \; x_2 = xyz, \; x_3 = xy^2, \; x_4 =
  y^2z, \; x_5 = xz^2, \; x_6 = yz^2\]
and a homogeneous ideal of $S$ in $\PN_6$ is given by all $2 \times 2$ minors of 
\[A = \left(\begin{array}{ccc}
          x_0 & x_1 & x_2 \\
          x_3 & x_2 & x_4 \\
          x_2 & x_5 & x_6
        \end{array}\right).\]

The ideal of $\PN_2 \times \PN_2$ embedded into $\PN_8$ by Segre is given
by the rank one locus of any general $3 \times 3$ matrix of linear
forms (if $s_1^p, s_2^p, s_3^p$ are homogeneous coordinates of the
$p$--th factor, write $x_{i,j} = s_i^1s_j^2$ for the homogeneous
coordinates of $\PN_8$). The rank one locus of such a general matrix over $\PN_7$ is
therefore isomorphic to a hyperplane section of $\PN_2 \times \PN_2$,
i.e., to $\PN(T_{\PN_2})$. The idea is hence to degenerate the general matrix into $A$.

The smoothing with general fiber $\PN_1 \times \PN_1 \times \PN_1$ is
obtained similarly: the ideal of
$S$ can as well be described as minors of the cube
\[\xymatrix@!0{ 
 & & x_1 \ar@{-}[rr]\ar@{-}'[d][dd]
  & & x_5 \ar@{-}[dd]\\
 & x_0 \ar@{-}[ur]\ar@{-}[rr]\ar@{-}[dd]
   & & x_2 \ar@{-}[ur]\ar@{-}[dd]\\
 A' = & & x_2 \ar@{-}'[r][rr]
  & & x_6\\
 & x_3 \ar@{-}[rr]\ar@{-}[ur]
  & & x_4 \ar@{-}[ur]}\]
What one has to do here is to compute all minors along side faces and diagonals. As above
we compare this with the {\em general} cube with linear entries over
$\PN_7$. The general cube cuts out $\PN_1 \times \PN_1
\times \PN_1$ (if $s_1^p, s_2^p$ are homogeneous coordinates of the
$p$--th factor, write $x_{i,j,k} = s_i^1s_j^2s_k^3$ for the
homogeneous coordinates of $\PN_7$).

Notice that the above description of the anticanonical embedding 
also shows that the general hyperplane section
through the vertex of $X$ has an elliptic Gorenstein
point, i.e., the vertex is a non--cDV point (as is any isolated canonical non--terminal point).
\end{example}


\end{document}